\documentclass[a4paper, oneside, reqno]{amsart}
\usepackage{amssymb,amsmath,mathtools}
\usepackage{graphicx}
\usepackage[labelformat=simple]{subcaption}
\usepackage{epstopdf}
\usepackage[usenames, dvipsnames]{color}
\calclayout
\usepackage{lscape}
\usepackage[colorlinks=true, urlcolor=blue, linkcolor=RedOrange, citecolor=TealBlue, pdfborder={0 0 0}]{hyperref}
\usepackage{enumerate}
\usepackage{multirow}
\begin{document}

\title[M. Nagy, S. R. Cowell and V. Beiu: Survey of Cubic Fibonacci Identities]{\textbf{Survey of Cubic Fibonacci Identities \\ When Cuboids Carry Weight}}
\author{\textbf{Mariana Nagy$^{0000-0001-6702-8049}$, Simon R. Cowell$^{0000-0001-6093-9755}$ and \\ Valeriu Beiu$^{0000-0001-8185-956X}$}}
\thanks{This work was supported in part by the European Union through the European Regional Development Fund under the Competitiveness Operational Program (\emph{BioCell-NanoART} = \emph{Novel Bio-inspired Cellular Nano-architectures}, POC-A1.1.4-E-2015 nr. 30/01.09.2016).}
%\address{Department of Mathematics and Computer Science\\ 
%``Aurel Vlaicu" University of Arad \\ 
%Str. Elena Dr\u{a}goi 2-4, RO-310330 Arad, Romania.}
%\email{mariana.nagy@uav.ro, simon.cowell@uav.ro, valeriu.beiu@uav.ro}
\date{\today}
\begin{abstract}
The aim of this paper is to present a comprehensive survey of cubic Fibonacci identities, trying to uncover as many as possible. From the outset, our rationale for a very careful search on an apparently obscure problem was not only a matter of mathematical curiosity, but also motivated by a quest for 3D Fibonacci spirals.

As we were not able to find any survey on the particular topic of cubic Fibonacci identities we decided to try to fill this void. We started by surveying many Fibonacci identities and recording cubic ones. Obviously, tracing all Fibonacci identities (for identifying a handful) is a daunting task. Checking several hundred we have realized that it is not always clear who the author is. The reason is that in many cases an identity was stated in one article (sometimes without a proof, e.g., as an open problem, or a conjecture) while later being proven in another article, or effectively rediscovered independently by other authors. However, we have done our best to present the identities chronologically. We have supplied our own proof for one identity, having tried, but failed, to find a published proof. For all the other identities, we either proved them on a computer, or else verified by hand their original published proofs. Somehow unexpectedly, our investigations have revealed only a rather small number of cubic Fibonacci identities, representing a tiny fraction of all published Fibonacci identities (most of which are linear or quadratic). Finally, out of these, only a handful of cubic Fibonacci identities are homogeneous. 
\end{abstract}
\maketitle

\section{Introduction}
Many scientists have been fascinated by numbers. They have relied on number theory for establishing properties and relations amongst digits, numbers and/or sequences of digits/numbers. Early on, such quests were mostly curiosity driven, while nowadays the interest is motivated by different computing problems, some of which are very practical. 

	One of the most famous sequences of numbers is the Fibonacci one. The fascination with this sequence started in the late 1800s, when \'{E}douard Lucas was calling it Lam\'{e}'s sequence. Since then, many have analyzed this sequence, discovering plenty of properties and identities, relating it to other sequences (e.g., Pell, Lucas, etc.), formulating quite a few generalizations, as well as finding various applications in sciences and engineering. Among all of these, a slightly less investigated aspect is represented by the relations between Fibonacci numbers and geometry. By far the most well-known result here is the 2D Fibonacci spiral, while relations to 3D seem to have been introduced for the first time by Johannes Kepler in 1611 \cite{Kepler_1611}. Therefore, a 3D version of the famous Fibonacci spiral was something to be expected. To our surprise, a search we have recently done \cite{Nagy_Cowell_and_Beiu_2018} has revealed that a 3D Fibonacci spiral is lacking, as none of the so-called 3D spirals are backed by a Fibonacci identity describing a recurrent tiling of 3D space. In fact, \cite{Nagy_Cowell_and_Beiu_2018} has shown that for distinguishing different tilings of the 3D space, a thorough understanding of cubic Fibonacci identities is needed. 

As we were not able to find any surveys of cubic Fibonacci identities we decided to start working on our own. That is why, in our journey through identities associated with the Fibonacci sequence, we have been looking for cubic ones. We have decided, whenever possible, to keep the original notations in the text. However, for uniformity in the summary tables towards the end of the paper, we have used a standard notation. 

The structure of the paper is as follows. In Section 2 we will introduce the Fibonacci numbers and present some of the early identities involving them. This will be followed in Section 3 by a very thorough survey of cubic Fibonacci identities. The cubic identities will be introduced chronologically and afterwards classified into: homogeneous, non-homogeneous, and generic, being reported in three different tables. For the identities reported in these tables, we first checked that generating their first few terms on a computer did not falsify them. We then verified them, either by a computer algebra system or else by checking the published proof by hand, and in one case, providing our own proof. The paper ends with conclusions and further directions of research.

\section{Fibonacci Numbers and Identities}

One of the most studied sequences of numbers was introduced by Leonardo of Pisa in his book \textit{Liber Abaci} \cite{Fibonacci_1202}. Leonardo of Pisa, or Leonardo Pisano, is better known today as Fibonacci \cite{Devlin_2011}. This is short for \textit{filius Bonacci} which means \textit{son of Bonacci}. He is regarded as one of the greatest and most creative minds of the Middle Age \cite{Devlin_2017}.

\textit{The Book of Calculation} is called \textit{Liber Abaci} \cite{Fibonacci_1202} and should be seen as ``liberating'' calculations, by showing how to use the Hindu/Arabic numerals and the place-valued decimal system. It is an extensive work based on the theoretical foundations of Euclidean mathematics, including arithmetic (number theory), algebra, and applied mathematics. \textit{Liber Abaci} was first published in 1202 (a second revised version appeared in 1228). A very thorough and long awaited English translation was done by Laurence E. Sigler and was published posthumously in 2002 \cite{Sigler_2002}, i.e., 800 years after \cite{Fibonacci_1202}.

\textit{Liber Abaci} starts by explaining the principle of the place-valued decimal system and demonstrating its use in arithmetic operations. Afterwards, the book refers to, and proposes mathematical solutions to, a plethora of practical problems pertaining to business and finance (e.g., profit margin, barter, money changing, conversion of weights and measures, partnerships, and interest), but also arts. By far the most well-known problem is the so-called \textit{Rabbit Problem}, proposed in Chapter Twelve (see page 404 in \cite{Sigler_2002}). The problem is to calculate how many pairs of rabbits (one male and one female) there will be at the end of one year, if starting with one pair of rabbits. In one month each pair bears another pair, and the newborn pair bears a pair starting from the second month. The answer is that after twelve months there will be 377 pairs of rabbits. The monthly evolution is described by the following sequence of numbers:
\begin{equation}
\label{eq: 1a}
1, 2, 3, 5, 8, 13, 21, 34, 55, 89, 144, 233, 377.
\end{equation}
This hypothetical breeding of rabbits introduces a sequence. This sequence of numbers (or their ratios) appears in nature (see, e.g. \cite{Mitchison_1977}): in pine cones, pineapples, artichokes, branching of trees, petals of flowers, arrangements of leaves, patterns of florets, and in the geometry of certain shells. It was the French mathematician \'{E}douard Lucas \cite{Lucas_1877} who first associated Fibonacci's name with this sequence.

It should be mentioned that these numbers seem to have been known for quite some time. In fact, historians of Indian mathematics have discovered many ideas that went unnoticed \cite{Singh_1936}. The main reasons were that they were: (i) written in Sanskrit; and (ii) not included in Western-style mathematical treatises, but embedded in \textit{s\={u}tras} in treatises on Sanskrit prosody (\textit{chandas}). This is the science of poetic meters and verse (arrangements of syllables) essential to the Vedas.

The oldest such treatise which has been preserved is Pi\.{n}gala's \textit{Chandaḥ\'{s}\={a}stra} \cite{Singh_1936}, \cite{Pingala_200BC}, which has 8 chapters containing 310 \textit{s\={u}tras}. It is a collection of aphorisms presenting a formal theory of poetic meters. The six topics of interest were: (i) listing all possible forms of a meter; (ii) recovering the form of a meter from its serial number; (iii) determining the serial number from a given form; (iv) calculating the number of forms having a given number of syllables; (v) calculating the total number of forms of a meter; (vi) determining the space needed to write all the forms of a meter. For solving these problems Pi\.{n}gala made use of a binary notation system, recursion, geometric series, sequences, partial sums, and binomial coefficients. These support claims that some \textit{s\={u}tras} from \cite{Pingala_200BC} suggest Pascal's triangle (known as \textit{meru prast\={a}ra}), while others link to Fibonacci numbers (called \textit{m\={a}tr\={a} meru}). It follows that, most probably, such numbers were known to Pi\.{n}gala (see \cite{Singh_1985} and \cite{Shah_2013}), but the rules forming them are hard to grasp from \cite{Pingala_200BC}. Such ideas were analyzed not only in relation to poetry, but also to music \cite{Hall_2008}, \cite{Wilson_1992}.

The number of patterns that are $m$ short syllables long (variations of \textit{m\={a}tr\={a}}-meters) is the ($m+1$)-th Fibonacci number. It is now accepted that the first algorithmic description appeared in the work of Virah\={a}\.{n}ka (between 600 and 800 AD) \cite{Shah_2013}, \cite{Velankar_1962}. Unfortunately, the original work was lost, being accessible only through the comments of Gop\={a}la on a palm leaf manuscript (c. 1133-1135) \cite{Singh_1985}, \cite{Shah_2013}, \cite{Velankar_1962}. This presents a recursion for the computation of the number of forms of \textit{m\={a}tr\={a}}-meters. A translation (with comments) of the algorithm reads: ``\textit{The count (of all possible forms of a m\={a}tr\={a}-meter) is obtained by adding (the counts of) permutations of the two previous (m\={a}tr\={a}-meters). This is the way to the count (of total permutations) of succeeding m\={a}tr\={a}s (m\={a}tr\={a}-meters)}'' \cite{Shah_2013}. Obviously, this gives the well-known recursion for what is now known as the Fibonacci sequence:
\begin{equation}
\label{eq: 1b}
F_{n+1}=F_n+F_{n-1}
\end{equation}
with the initial conditions $F_0=0$ and $F_1=1$ (see \cite{OEIS_A000045}).

A first historical perspective written by Brooke \cite{Brooke_1964} on Fibonacci numbers covered till 1900. It mentioned that Johann Kepler rediscovered the recurrence in 1611 \cite{Kepler_1611}, and that Jacques Binet derived a closed analytical form in 1843 \cite{Binet_1843}. In fact, the analytical form was proven by Abraham de Moivre one century earlier \cite{deMoivre_1730} (see Chp. 2: \textit{De natura Serierum recurrentium}, pp. 26‒42). The sequence \eqref{eq: 1a} is also mentioned by Simon Stevin in 1634 (\cite{Stevin_1634} on page 170), and by Robert Simson in 1753 (\cite{Simson_1753} on page 368), who refers to Albert Girard's commentary on Stevin's work \cite{Stevin_1634}.

The growing fascination with these numbers started properly only in the late 1800's, after \'{E}douard Lucas published \cite{Lucas_1877}. Since then many scientists have focused on this sequence and an impressive number of properties and identities have been discovered. The sequence was studied also in relation with other known sequences of numbers (e.g., Pell \cite{OEIS_A000129} and Lucas numbers \cite{OEIS_A000032}), and various generalizations have been defined.

Before going further, we shall enumerate only a few Fibonacci identities. The earliest ones are from 1680 and are due to Cassini \cite{Cassini_1680} (presented as text, without notation) %, on page 201
\begin{equation}
\label{eq: 2a}
F_{n-1} F_{n+1}-F_n^2=(-1)^n
\end{equation}
\begin{equation}
\label{eq: 2b}
F_{n+2} F_{n-1}-F_n F_{n+1}=(-1)^n
\end{equation}
\begin{equation}
\label{eq: 2c}
F_{n-2} F_{n+2}-F_n^2=(-1)^{n+1}.
\end{equation}
These were generalized by Catalan \cite{Catalan_1886} in 1879 (and mentioned in \cite{Dickson_1919} on page 402):
\begin{equation}
\label{eq: 3}
U_n^2-U_{n-p} U_{n+p}=(-1)^{n+1-p} U_{p-1}^2,
\end{equation}
by d'Ocagne \cite{dOcagne_1886} in 1885:
\begin{equation}
\label{eq: 4}
u_p u_i-u_{p+1} u_{i-1}=(-1)^{i+1} u_{p-i+1}
\end{equation}
and by Tagiuri \cite{Tagiuri_1901} in 1901:
\begin{equation}
\label{eq: 5}
u_s u_n-u_{s+k} u_{n-k}=(-1)^{n-k} u_k u_{s-n+k}.
\end{equation}

Lucas added a few of his own in 1876 \cite{Lucas_1876} and 1877 \cite{Lucas_1877} (see also \cite{Dickson_1919}):
\begin{equation}
\label{eq: 6a}
U_{n-1}^2+U_n^2=U_{2n}
\end{equation}
\begin{equation}
\label{eq: 6b}
F_{2n+1}=F_{n+1}^2+F_n^2
\end{equation}
\begin{equation}
\label{eq: 6c}
F_{2n}=F_{n+1}^2-F_{n-1}^2
\end{equation}
\begin{equation}
\label{eq: 6d}
F_{n+3}^2=2F_{n+2}^2+2F_{n+1}^2-F_n^2
\end{equation}
\begin{equation}
\label{eq: 6e}
\sum_{k=0}^n F_k =F_{n+2}-1
\end{equation}
\begin{equation}
\label{eq: 6f}
\sum_{k=0}^n F_k^2=F_n F_{n+1}
\end{equation}
and a few more already known and appreciated. As an example, \eqref{eq: 6f} supports the construction of the well-known two dimensional Fibonacci spiral.

Several authors identified classes of identities and proposed general formulas for them. Such identities are given by Lucas \cite{Lucas_1877} as symbolic formulas where, after expansion, exponents are replaced by subscripts:
\begin{equation}
\label{eq: 7a}
U^{n+2p}=U^n (U+1)^p
\end{equation}
\begin{equation}
\label{eq: 7b}
U^{n-p}=U^n (U-1)^p.
\end{equation}
For different values of $p$, we obtain valid linear identities; in particular, setting $p=1$ in \eqref{eq: 7a} or \eqref{eq: 7b} yields \eqref{eq: 1b}.

Another formula, published in 1965 by Halton \cite{Halton_1965}:
\begin{equation}
\label{eq: 8a}
F_m^k F_n=(-1)^{kr} \sum_{h=0}^k C_k^h (-1)^h F_r^h F_{r+m}^{k-h} F_{n+kr+hm}
\end{equation}
is general as it establishes a class of identities, being valid for any $(k,r,m,n)$ quadruplet of non-negative integers, where $C_k^h$ is the number of $h$-combinations from a set of $k$ elements. In fact, equation \eqref{eq: 8a} was reported earlier in a different form in \cite{Vinson_1961}. As an example, for $(2,1,1,n)$ we obtain: 
\begin{equation}
\label{eq: 8b}
F_n=F_{n+2}-2F_{n+3}+F_{n+4}.
\end{equation}

\section{Cubic Fibonacci Identities}

Chronologically, the \textbf{first cubic Fibonacci identity} we have been able to find is the one mentioned in 1919 by Dickson in his \textit{History of the Theory of Numbers} \cite{Dickson_1919}. In Chapter XVII ``Recurring Series: Lucas' $u_n, v_n$'' (volume I) the following identity appears on page 395:
\begin{equation}
\label{eq: 9a}
U_n^3+U_{n+1}^3-U_{n-1}^3=U_{3n+2}.
\end{equation}
The notation $U_n$ used by Dickson starts with $U_1=1$, $U_2=2$ as initial conditions, leading to $1, 2, 3, 5, 8, \ldots$. Dickson mentioned that \eqref{eq: 9a} was stated without proof by Lucas in \cite{Lucas_1876}. In fact, in \cite{Lucas_1877} on page 11, Lucas shows that $u_n^3+u_{n+1}^3-u_{n-1}^3=u_{3n}$, with $u_{n+2}=u_{n+1}+u_n$ and initial conditions $u_0=0$, $u_1=1$, which gives $0, 1, 1, 2, 3, 5, 8, \ldots$. We were not expecting that this would be the one and only cubic Fibonacci identity \cite{OEIS_A056570} mentioned in \cite{Dickson_1919}. 

Over time, this cubic Fibonacci identity has resurfaced more than once being mentioned by, e.g., Vorob'ev \cite{Vorobev_1961} (using Lucas' notation) and Subba Rao \cite{Subba_Rao_1953} (using Dickson's notation). Later, it was mentioned briefly in \cite{Brousseau_1972}, which pointed to a note signed just G. that ``\textit{quoted the cubic relation given by Dickson}.'' Brousseau states in \cite{Brousseau_1972} that G. is ``\textit{evidently J. Ginsburg himself}'' (the Editor of \textit{Scripta Mathematica} \cite{Ginsburg_1953}). Still, in \cite{Brousseau_1972} Brousseau uses the notation:
\begin{equation}
\label{eq: 9b}
F_{n+1}^3+F_n^3-F_{n-1}^3=F_{3n}
\end{equation}
which is the one we shall retain in this paper (also being used by OEIS, see \cite{OEIS_A000045} and \cite{OEIS_A056570}).

More recently, the book of Benjamin \& Quinn \cite{Benjamin_and_Quinn_2003} includes this cubic Fibonacci identity as one of the $246$ identities covered in their book. It appears as Identity $232$ in Chapter 9 ``Advanced Fibonacci \& Lucas Identities'' being written as $f_{n-1}^3+f_n^3-f_{n-2}^3=f_{3n-1}$ (where $f_n= F_{n+1}$) and proven on page 126.

Lucas also presented \textbf{a second cubic Fibonacci identity} in \cite{Lucas_1877} on pages 16‒17:
\begin{equation}
\label{eq: 10}
5U_n^3=U_{3n}-3(-1)^n U_n.
\end{equation}

This identity is mentioned by Halton \cite{Halton_1965} as a special case deduced from \eqref{eq: 8a} and, later, it is used by Melham \cite{Melham_2008} while stating several conjectures on the sums of powers of Fibonacci numbers.

\textbf{A third cubic Fibonacci identity} appeared in a book published in 1951 in Russian by Vorob'ev \cite{Vorobev_1961}. It is \textit{the first cubic Fibonacci identity involving a ``sum of Fibonacci cubes''} (i.e., $\sum F_i^3$ ) and was written using Lucas's notation:
\begin{equation}
\label{eq: 11}
u_1^3+u_2^3+ \cdots +u_n^3=\frac{1}{10} \left [ u_{3n+2}+(-1)^{n+1} 6u_{n-1}+5 \right ].
\end{equation}
The book was translated to English ten years later. The same identity was proved by Clary \& Hemenway \cite{Clary_and_Hemenway_1993} (using the $F_n$ notation), while later being stated (without proof) by Benjamin \& Quinn \cite{Benjamin_and_Quinn_2003} (on page 14).

\textbf{The fourth cubic Fibonacci identity} seems to be the one introduced by Block in a short note \cite{Block_1953}. It is the first cubic Fibonacci identity which involves a ``sum of Fibonacci cuboids'' (e.g., $\sum F_i F_j F_k$ and/or $\sum F_i^2 F_j$) while also being \textit{the first homogeneous one}: 
\begin{equation}
\label{eq: 12}
\sum_{k=1}^n F_k^2  F_{k+1}=\frac{1}{2}F_n F_{n+1} F_{n+2}.
\end{equation}
Zeitlin wrote it using $f_n$ and proved it by mathematical induction in \cite{Zeitlin_1963} (as identity XIX). Harris proved it once again, using the method of summation by parts \cite{Harris_1965}, while Pond proved a generalized version in \cite{Pond_1968}. Afterwards, it was mentioned by Clary \& Hemenway \cite{Clary_and_Hemenway_1993} and Melham \cite{Melham_1999}.

An identity similar to \eqref{eq: 9b} was introduced in December 1953 \cite{Ginsburg_1953}:
\begin{equation}
\label{eq: 13}
F_{n+2}^3-3F_n^3+F_{n-2}^3=3F_{3n}
\end{equation}
being proved by Halton \cite{Halton_1965} as a special case of \eqref{eq: 8a}, and later stated by Brousseau in \cite{Brousseau_1972}.

Simultaneously, Subba Rao (using Dickson's notation \cite{Dickson_1919}) published a gold mine of cubic Fibonacci identities \cite{Subba_Rao_1953}. Besides mentioning one of the early Lucas identities, he adds ten new cubic Fibonacci identities. These can be classified into:
\begin{itemize}
\item
two classical identities
\begin{equation}
\label{eq: 14}
U_m U_{m+1} U_{m+2}=U_{m+1}^3+(-1)^{m+1} U_{m+1}
\end{equation}
\begin{equation}
\label{eq: 15}
U_{m+1}^3=U_m^3+U_{m-1}^3+3U_{m-1} U_m U_{m+1};
\end{equation}
\item
three including sums of cubes
\begin{equation}
\label{eq: 16}
\sum_{k=0}^n U_{2k}^3 =\frac{1}{4} \left (U_{2n+1}^3+3U_{2n+1} \right )
\end{equation}
\begin{equation}
\label{eq: 17}
\sum_{k=1}^n U_{2k-1}^3 = \frac{1}{4} \left (U_{2n}^3-3U_{2n}+2 \right )
\end{equation}
\begin{equation}
\label{eq: 18}
\sum_{k=1}^n U_k^3 = \frac{1}{4} [U_{n+1}^3+U_n^3+(-1)^n 3U_{n-1}-2];
\end{equation}
\item
five including sums of cuboids
\begin{equation}
\label{eq: 19}
\sum_{k=1}^n U_{2k-1} U_{2k} U_{2k+1} = \frac{1}{4} (U_{2n+1}^3-U_{2n+1} )
\end{equation}
\begin{equation}
\label{eq: 20}
\sum_{k=1}^n U_{2k-2} U_{2k-1} U_{2k} = \frac{1}{4} (U_{2n}^3+U_{2n}-2)
\end{equation}
\begin{equation}
\label{eq: 21}
\sum_{k=1}^n U_{2k-2} U_{2k} U_{2k+2} = \frac{1}{4} (U_{2n+1}^3+7U_{2n+1}-8)
\end{equation}
\begin{equation}
\label{eq: 22}
\sum_{k=1}^n U_{2k-1} U_{2k+1} U_{2k+3} = \frac{1}{4} (U_{2n+2}^3-7U_{2n+2}+6)
\end{equation}
\begin{equation}
\label{eq: 23}
\sum_{k=1}^n U_{k-1} U_{k+1} U_{k+3} = \frac{1}{10} [U_{3n+7}+(-1)^{n+1} 16U_n-5].
\end{equation}
\end{itemize}

Later, in his search for the properties of certain generalized functions, Horadam \cite{Horadam_1965} gave an algebraic proof of \eqref{eq: 18}. In \cite{Koshy_2001} Koshy mentioned \eqref{eq: 15}, attributing it to Carlitz (but the reference seems to be missing). Clary \& Hemenway \cite{Clary_and_Hemenway_1993} mentioned \eqref{eq: 17}, and used it as an intermediary step for proving an identity involving both Fibonacci and Lucas numbers.

 \textbf{The first two alternating sums} seem to be the ones introduced by Zeitlin in 1963 \cite{Zeitlin_1963}, both for cuboids (identity XX) and for cubes (identity XXI):
\begin{equation}
\label{eq: 24}
2 \sum_{k=0}^n (-1)^k f_k f_{k+1}^2 = (-1)^n f_n f_{n+1} f_{n+2}
\end{equation}
\begin{equation}
\label{eq: 25}
2 \sum_{k=0}^n (-1)^k f_{k+1}^3 = (-1)^n [f_{n+1}^2 f_{n+4} -f_n f_{n+2} f_{n+3} ]-1 = f_{n+4}+(-1)^n f_n f_{n+2}^2 - 1.
\end{equation}

Hogatt \cite{Hoggatt_1963_1}, \cite{Hoggatt_1963_2} recorded a general non-homogeneous identity:
\begin{equation}
\label{eq: 26}
F_{r+1} F_{s+1} F_{t+1} + F_r F_s F_t + F_{r-1} F_{s-1} F_{t-1} = F_{r+s+t}
\end{equation}
(proposed by I. D. Ruggles), as well as a homogeneous one:
\begin{equation}
\label{eq: 27}
F_{n+1}^3 = 3F_n^3+6F_{n-1}^3-3F_{n-2}^3-F_{n-3}^3
\end{equation}
(proposed by Zeitlin \& Parker). Only a few months later Brown \cite{Hoggatt_1963_3} proved \eqref{eq: 26}, while later Johnson, using a matrix representation, presented an alternate proof \cite{Johnson_2009}. A proof of \eqref{eq: 27} appeared in \cite{Hoggatt_1963_4}, while it was also presented by Brousseau \cite{Brousseau_1968} along with a few other cubic identities. Recently \eqref{eq: 27} links to \cite{Khomovsky_2018}, being a special case of classes of Fibonacci identities (see later).

Another general identity was proposed by Halton \cite{Halton_1965} (see \eqref{eq: 8a}), and was supplemented by Gould \cite{Gould_1964}:
\begin{equation}
\label{eq: 28}
F_{n+a} F_{n+b} F_{n+c}-F_n F_{n+a} F_{n+b+c}+F_n F_{n+b} F_{n+c+a}-F_n F_{n+c} F_{n+a+b}=(-1)^n (F_a F_b F_{n+c}-F_c F_a F_{n+b}+F_b F_c F_{n+a} )
\end{equation}
(see also \cite{Fairgreave_and_Gould_2005}). This can be used to generate many cubic homogeneous Fibonacci identities.

To this growing list, having Subba Rao \cite{Subba_Rao_1953} and Zeitlin \cite{Zeitlin_1963} as major contributors, Harris \cite{Harris_1965} added two more cubic Fibonacci identities, one including sums of cuboids:
\begin{equation}
\label{eq: 29}
\sum_{k=0}^n f_k^2 f_{k+2} = \frac{1}{2} [f_{n+3} f_{n+1} f_n-(-1)^n f_{n-1}+1],
\end{equation}
and another one including sums of cubes:
\begin{equation}
\label{eq: 30}
\sum_{k=0}^n f_k^3 = \frac{1}{2} [f_{n+1}^2 f_n+(-1)^n f_{n-1}+1].
\end{equation}
This one was later mentioned by Brousseau in \cite{Brousseau_1972}.

As if not to be outdone, Pond detailed and proved quite a few cubic identities for $H_n$ \cite{Pond_1968}. $H_n$ is a generalized Fibonacci sequence $H_n=H_{n-1}+H_{n-2}$ having arbitrary initial conditions. For $H_0=0$ and $H_1=1$ one gets the classical Fibonacci sequence, and \eqref{eq: 12} resulted as a byproduct of the cubic identities proved in \cite{Pond_1968}. We also mention here a homogeneous cubic Fibonacci identity from \cite{Pond_1968}:
%\begin{equation}
%\label{eq: 31}
%\sum H_{n-1} H_n H_{n+1}+ \sum H_n^3 = H_n H_{n+1}^2
%\end{equation}
\begin{equation}
\label{eq: 31}
\sum_{k=1}^n H_{k-1} H_k H_{k+1} + \sum_{k=0}^n H_k^3 = H_n H_{n+1}^2
\end{equation}
but there are several variations in \cite{Pond_1968} at intermediate steps of the proofs. For both \eqref{eq: 30} and \eqref{eq: 31} combinatorial proofs were given by Benjamin, Carnes \& Cloitre \cite{Benjamin_Carnes_and_Cloitre_2009}. Identity \eqref{eq: 31} appeared also in \cite{Lang_2012}, along with a few other cubic Fibonacci identities we shall mention later.

Taking a straightforward approach, Brousseau introduced several identities by cubing well-known linear ones \cite{Brousseau_1968}:
\begin{equation}
\label{eq: 32}
F_{n+1}^3=F_n^3+3F_n^2 F_{n-1}+3F_n F_{n-1}^2+F_{n-1}^3
\end{equation}
\begin{equation}
\label{eq: 33}
F_{n-2}^3=F_n^3-3F_n^2 F_{n-1}+3F_n F_{n-1}^2-F_{n-1}^3
\end{equation}
\begin{equation}
\label{eq: 34}
F_{n-3}^3=-F_n^3+6F_n^2 F_{n-1}-12F_n F_{n-1}^2+8F_{n-1}^3
\end{equation}
and the list continues by rediscovering \eqref{eq: 27} and several others. Obviously, as cubes of linear identities, all these cubic Fibonacci identities are homogeneous.

Our quest for cubic identities led us to several more: 
\begin{equation}
\label{eq: 35}
F_{3n}=F_{n+2}^3-F_{n-1}^3-3F_n F_{n+1} F_{n+2}
\end{equation}
\begin{equation}
\label{eq: 36}
F_n F_{n+3}^2-F_{n+2}^3=(-1)^{n+1} F_{n+1}
\end{equation}
\begin{equation}
\label{eq: 37}
F_{n+3} F_n^2-F_{n+1}^3=(-1)^{n+1} F_{n+2}.
\end{equation}
Equation \eqref{eq: 35} was proposed by Padilla being mentioned in \cite{Hillman_1970} along with \eqref{eq: 15}, while \eqref{eq: 36} and \eqref{eq: 37} were introduced by Hogatt \& Bergum \cite{Hoggatt_and_Bergum_1977}.

New Fibonacci identities were stated and/or proved by Clary \& Hemenway \cite{Clary_and_Hemenway_1993}, including two cubic ones: 
\begin{equation}
\label{eq: 38}
\sum_{k=1}^n F_k^3 =\frac{1}{10} [F_{3n+2}-(-1)^n 6F_{n-1}+5]
\end{equation}
\begin{equation}
\label{eq: 39}
\sum_{k=1}^n F_{2k}^3 =\frac{1}{4} (F_{2n+1}-1)^2 (F_{2n+1}+2).
\end{equation}

After 2000, the first significant contribution is \cite{Koshy_2001}, where Koshy presents a collection of Fibonacci and Lucas sequences, the real \textit{``twin shining stars''} of integer sequences. Besides proving many properties and showing applications in various fields such as art, biology, chemistry, geometry, etc., the book lists a large number of both proved and to-be-proved Fibonacci identities. Among these identities, one can easily recognize \eqref{eq: 14}-\eqref{eq: 23} and \eqref{eq: 27}.

Two years later, Melham \cite{Melham_2003} dedicates a whole paper to stating and proving the following cubic identity:
\begin{equation}
\label{eq: 40}
F_{n+1} F_{n+2} F_{n+6}-F_{n+3}^3=(-1)^n F_n,
\end{equation}
while in 2004, Benjamin, Carnes \& Cloitre published \cite{Benjamin_Carnes_and_Cloitre_2009}, that mentions or proves many Fibonacci, Lucas and related generalized identities, including \eqref{eq: 27} and \eqref{eq: 38}.

Fairgrieve \& Gould \cite{Fairgreave_and_Gould_2005} ran extensive computer searches and discovered a companion formula to \eqref{eq: 40}:
\begin{equation}
\label{eq: 41}
F_n F_{n+4} F_{n+5}-F_{n+3}^3=(-1)^{n+1} F_{n+6},
\end{equation}
along with another pair of dual formulas, for which they also give combinatorial proofs:
\begin{equation}
\label{eq: 42}
F_{n-2} F_{n+1}^2-F_n^3=(-1)^{n-1} F_{n-1}
\end{equation}
\begin{equation}
\label{eq: 43}
F_{n+2} F_{n-1}^2-F_n^3=(-1)^n F_{n+1}.
\end{equation}

As an interesting alternative to traditional algebraic or combinatorial proofs, Ollerton \cite{Ollerton_2006} uses a 3D graphical approach which has clear educational advantages (for 2D examples the reader should see \cite{Brousseau_1972}). Ollerton explains his method and uses it to prove \eqref{eq: 27} and proposes:
\begin{equation}
\label{eq: 44}
F_{l+1} F_{m+1} F_{n+1}=3F_l F_m F_n+6F_{l-1} F_{m-1} F_{n-1}-3F_{l-2} F_{m-2} F_{n-2}-F_{l-3} F_{m-3} F_{n-3}.
\end{equation}

In 2008, Benjamin, Eustis \& Plott \cite{Benjamin_Eustis_and_Plott_2008} provide combinatorial proofs for a dozen identities stated without proof in \cite{Benjamin_and_Quinn_2003}, while adding a general formula of their own:
\begin{equation}
\label{eq: 45}
\sum_{k=0}^n C_n^k f_{qk-1+p} (f_{q-3} )^{n-k}=(f_{q-1} )^n f_{2n-1+p}.
\end{equation}
Taking $p=0$ and $n=3$ we get another cubic identity:
\begin{equation}
\label{eq: 46}
8F_n^3=3F_n F_{n-2}^2+3F_{2n} F_{n-2}+F_{3n}.
\end{equation}

Chen \& Chen’s approach in \cite{Chen_and_Chen_2008} and \cite{Chen_2010} relies on proving that for any positive integers $k$ and $n$ one can find a homogeneous identity consisting of a linear combination of powers of the Fibonacci numbers:
\begin{equation}
\label{eq: 47}
F_{n+k+1}^k=\sum_{i=0}^k a_i (k)F_{n+k-i}^k
\end{equation}
where $a_i(k)$ are integers. They also show how to determine $a_i(k)$. For $k=3$ one gets \eqref{eq: 27}.

More than a century after it was introduced by Lucas \cite{Lucas_1877}, and several decades after it was published by Subba Rao \cite{Subba_Rao_1953}, \eqref{eq: 10} was used by Melham \cite{Melham_2008} as a starting point for a set of conjectures involving both Lucas and Fibonacci numbers. Proving those conjectures led Melham to also rediscover \eqref{eq: 27}. Later, Sun, Xie and Yang \cite{Sun_Xie_and_Yang_2016} mentioned \eqref{eq: 10}, and also proved one of Melham's conjectures.

In spite of such long standing efforts, in the last decade further scrutiny of the Fibonacci sequence went unabated. Thus, Melham \cite{Melham_2011} introduced several generic identities and, among others, established:
\begin{equation}
\label{eq: 48}
F_{n-3} F_{n+1}^2-F_{n-2}^2 F_{n+2}=4(-1)^n F_n
\end{equation}
along with its companion formula.

Azarian \cite{Azarian_2012_1}, \cite{Azarian_2012_2} proved the following two cubic identities by using binomial sums:
\begin{equation}
\label{eq: 49}
F_{n+1} [5F_{n+1}^2-3(-1)^n ]=F_{n+3}^3-F_n^3-3F_{n+1} F_{n+2} F_{n+3}=F_{n+1}^3+ F_{n+2}^3-F_n^3
\end{equation}
\begin{equation}
\label{eq: 50}
F_{3n}+4F_{3n+3}=  \frac{1}{3} ( F_{n+4}^3+F_n^3-3F_{n+2}^3 ).
\end{equation}

Almost simultaneously, Yazlik, Yilmaz \& Taskara \cite{Yazlik_Yilmaz_and_Taskara_2012} proved a new identity involving powers of generalized Fibonacci numbers $F_{k,n}$:
\begin{equation}
\label{eq: 51}
F_{k,n+4}^3=(k^3+2k) F_{k,n+3}^3+(k^4+3k^2+2) F_{k,n+2}^3-(k^3+2k) F_{k,n+1}^3-F_{k,n}^3,
\end{equation}
where the $F_{k,n}$ are defined by $F_{k,n+2}=k F_{k,n+1}+F_{k,n}$ with initial conditions $F_{k,0}=0$ and $F_{k,1}=1$.
For $k=1$ this reduces to \eqref{eq: 27}.

Lang \cite{Lang_2012}, while using partial summation (as an intermediate step) of cubes of Fibonacci numbers, mentions \eqref{eq: 31} while also proving two other homogeneous cubic identities:
\begin{equation}
\label{eq: 52}
\sum_{k=0}^n F_{k+2} F_{k+1} F_k =-\sum_{k=0}^n F_k^3 +F_{n+3} F_{n+1} F_n
\end{equation}
\begin{equation}
\label{eq: 53}
\sum_{k=0}^n F_{k+1} F_k^2 =\sum_{k=0}^n F_k^3 +\sum_{k=0}^{n-1} F_{k+1}^2 F_k.
\end{equation}

Lately, the focus has been on generalization of (known) identities and methods for computing terms of generalized Fibonacci formulas, but as far as we know, no fresh cubic identity was discovered. Still, it is worth mentioning Lang \& Lang \cite{Lang_and_Lang_2013} giving alternate proofs for several known identities (including \eqref{eq: 3}, \eqref{eq: 4} and \eqref{eq: 40}), and using recurrence relations to obtain:
\begin{equation}
\label{eq: 54}
F_{3n+3}=F_n F_{n+3}^2+F_{n-1} F_{n+2}^2- F_{n-2} F_{n+1}^2.
\end{equation}

Recently, Khomovsky \cite{Khomovsky_2018} presented a general method for obtaining recurrent identities. The Fibonacci sequence is used as an example. Besides obtaining \eqref{eq: 3}, a general identity is:
\begin{equation}
\label{eq: 55}
F_{2l} (F_{k+m}^2-F_{k-m}^2)=F_{2m} (F_{k+l}^2-F_{k-l}^2).
\end{equation}
For $l=k$ and $m=k+1$ this gives:
\begin{equation}
\label{eq: 56}
F_n F_{n+1}^2-F_n^2 F_{n+2}= (-1)^n F_n.
\end{equation}
Khomovsky also proves a sequence of identities involving Fibonacci and Lucas numbers, out of which classes of cubic Fibonacci identities can be obtained. Special cases are:
\begin{equation}
\label{eq: 57}
F_{n+2}^3+F_{n-2}^3=3(F_{n+1}^3-F_{n-1}^3 )+6F_n^3
\end{equation}
\begin{equation}
\label{eq: 58}
F_{n+3}^3-F_{n-3}^3=16(F_{n+1}^3-F_{n-1}^3 )+12F_n^3
\end{equation}
\begin{equation}
\label{eq: 59}
16F_{3n}=F_{n+3}^3+4F_n^3-F_{n-3}^3
\end{equation}
\begin{equation}
\label{eq: 60}
63F_{3n}=F_{n+4}^3-7F_n^3+F_{n-4}^3.
\end{equation}

Other cubic identities can be found in online resources such as Wikipedia \cite{Wikipedia} and Wolfram \cite{Chandra_and_Weisstein}. Chandra and Weisstein \cite{Chandra_and_Weisstein} can be credited with quite a selection of these, including:
\begin{equation}
\label{eq: 61}
F_{3n}=2F_n^3+3F_{n-1} F_n F_{n+1}
\end{equation}
\begin{equation}
\label{eq: 62}
F_{3n+1}=F_{n+1}^3+3F_{n+1} F_n^2-F_n^3
\end{equation}
\begin{equation}
\label{eq: 63}
F_{3n+2}=F_{n+1}^3+3F_n F_{n+1}^2+F_n^3,
\end{equation}
as well as a general one:
\begin{equation}
\label{eq: 64}
F_{kn+m}=\sum_{i=0}^k C_k^i F_{m-i} F_n^i F_{n+1}^{k-i},
\end{equation}
which they attribute to Mihailovs (personal communication). Taking $k=3$ and $m=3$ this gives:
\begin{equation}
\label{eq: 65}
F_{3n+3}=2F_{n+1}^3+3F_n F_{n+1}^2+3F_n^2 F_{n+1}.
\end{equation}

At this stage, all the cubic identities we have identified were rewritten according to the OEIS notation. Afterwards, we used Excel to ensure that each identity could not be falsified by checking its first several terms. All the identities having passed this test, we then used either Mathematica (which in some cases took over three hours), or else checked the published proofs to verify each identity in general. We also proved \eqref{eq: 64} in \cite{Cowell_Nagy_and_Beiu_2018}, as we were not able to find a published version. Finally, all the cubic identities were divided into three tables:
\begin{itemize}
\item
a first table with homogeneous cubic identities; 
\item
a second table with non-homogeneous cubic identities; 
\item
a third table listing general formulas that, in particular cases, can lead to cubic Fibonacci identities. 
\end{itemize}

\begin{table}%[!h]
%\label{tab: I}
% Previously the table labels came immediately after the \begin{table} commands, and the \refs to the tables were messed up. I discovered on stackexchange that \label must come after caption. Fixing this solved the problem.
\centering
\caption{Homogeneous Cubic Fibonacci Identities (In Chronological Order)}
\label{tab: I}
\resizebox{\textwidth}{!}{
\begin{tabular}{|c|c|c|l|}
\hline
\hline
Year & Identity & Condition & Author(s) \\
\hline
\hline
\multirow{9}{*}{1953} &
\multirow{6}{*}{$2\sum_{k=0}^n F_k^2 F_{k+1}=F_n F_{n+1} F_{n+2}$} &
\multirow{6}{*}{$n\geq0$	} &
Block \\
&&& \hspace{1cm} \emph{1963 \hspace{2mm} Zeitlin} \\
&&& \hspace{1cm} \emph{1965 \hspace{2mm} Harris} \\
&&& \hspace{1cm} \emph{1968 \hspace{2mm} Pond} \\
&&& \hspace{1cm} \emph{1993 \hspace{2mm} Clary \& Hemenway} \\
&&& \hspace{1cm} \emph{1999 \hspace{2mm} Melham} \\
\cline{2-4} &
\multirow{3}{*}{$F_{n+1}^3=F_n^3+F_{n-1}^3+3F_{n-1} F_n F_{n+1}$} &
\multirow{3}{*}{$n\geq1$} &
Subba Rao \\
&&& \hspace{1cm} \emph{1967 \hspace{2mm} Carlitz} \\
&&& \hspace{1cm} \emph{1970 \hspace{2mm} Hillman (Padilla)} \\
\hline
\multirow{11}{*}{1963}
& $2\sum_{k=0}^n (-1)^k F_k F_{k+1}^2=(-1)^n F_n F_{n+1} F_{n+2}$ &
\multirow{2}{*}{$n\geq0$} &
\multirow{2}{*}{Zeitlin} \\
\cline{2-2} &
% Use \hspace on either side of this equation, which is the widest in the table, to give it a little more room left-right, so as to better imitate Valeriu's original table:
\hspace{5mm}
$2\sum_{k=0}^n (-1)^k F_{k+1}^3=(-1)^n [F_{n+1}^2 F_{n+4}-F_n F_{n+2} F_{n+3}]-1$
\hspace{5mm} && \\
\cline{2-4} &
\multirow{9}{*}{$F_{n+4}^3+3F_{n+1}^3+F_n^3=3F_{n+3}^3+6F_{n+2}^3$} &
\multirow{9}{*}{$n\geq0$} &
Zeitlin \& Parker \\
&&& \hspace{1cm} \emph{1968 \hspace{2mm} Brousseau} \\
&&& \hspace{1cm} \emph{2001 \hspace{2mm} Koshy} \\
&&& \hspace{1cm} \emph{2003 \hspace{2mm} Benjamin \& Quinn} \\
&&& \hspace{1cm} \emph{2006 \hspace{2mm} Ollerton} \\
&&& \hspace{1cm} \emph{2008 \hspace{2mm} Chen \& Chen} \\
&&& \hspace{1cm} \emph{2010 \hspace{2mm} Chen} \\
&&& \hspace{1cm} \emph{2012 \hspace{2mm} Yazlik, Yilmaz \& Taskara} \\
&&& \hspace{1cm} \emph{2018 \hspace{2mm} Khomovsky} \\
\hline
\multirow{6}{*}{1968} &
\multirow{3}{*}{$\sum_{k=0}^n F_k^3+\sum_{k=1}^n F_{k-1} F_k F_{k+1}=F_n F_{n+1}^2$} &
\multirow{3}{*}{$n\geq1$} &
Pond \\
&&& \hspace{1cm} \emph{2004 \hspace{2mm} Benjamin, Carnes \& Cloitre} \\
&&& \hspace{1cm} \emph{2012 \hspace{2mm} Lang} \\
\cline{2-4}
& $F_{n+1}^3=F_n^3+3F_n^2 F_{n-1}+3F_n F_{n-1}^2+F_{n-1}^3$ &
$n\geq1$ &
\multirow{3}{*}{Brousseau} \\
\cline{2-3}
& $F_{n-2}^3=F_n^3-3F_n^2 F_{n-1}+3F_n F_{n-1}^2-F_{n-1}^3$ &
$n\geq2$ & \\
\cline{2-3}
& $F_{n-3}^3=-F_n^3+6F_n^2 F_{n-1}-12F_n F_{n-1}^2+8F_{n-1}^3$ &
$n\geq3$ & \\
\hline
2004
& $F_{n+1}F_{n+2}^2=F_{n+1}^3+\sum_{k=0}^n (F_k^3+F_k F_{k+1}F_{k+2})$ &
$n\geq0$ &
Benjamin, Carnes \& Cloitre \\
\hline
\multirow{2}{*}{2012}
& $\sum_{k=0}^n F_{k+2}F_{k+1}F_k +\sum_{k=0}^n F_k^3=F_{n+3}F_{n+1}F_n$ &
$n\geq0$ &
\multirow{2}{*}{Lang} \\
\cline{2-3}
& $\sum_{k=0}^n F_{k+1}F_k^2=\sum_{k=0}^n F_k^3+\sum_{k=0}^{n-1} F_{k+1}^2 F_k$ &
$n\geq1$ & \\
\hline
2018
& $F_{n+3}^3-F_{n-3}^3=16(F_{n+1}^3-F_{n-1}^3)+12F_n^3$ &
$n\geq3$ &
Khomovsky \\
\hline
\end{tabular}
}
\end{table}

% This is just a comment to separate the two tables, so that hopefully I will be less likely to get confused as to which table I'm editing !
% This is just a comment to separate the two tables, so that hopefully I will be less likely to get confused as to which table I'm editing !
% This is just a comment to separate the two tables, so that hopefully I will be less likely to get confused as to which table I'm editing !
% This is just a comment to separate the two tables, so that hopefully I will be less likely to get confused as to which table I'm editing !

%\newpage

\begin{table}%[!h]
%\label{tab: II}
\centering
\caption{Non-Homogeneous Cubic Fibonacci Identities (In Chronological Order)}
\label{tab: II}
\resizebox{\textwidth}{!}{
\begin{tabular}{|c|c|c|l|}
\hline
\hline
Year & Identity & Condition & Author(s) \\
\hline
\hline
\multirow{7}{*}{1876} &
\multirow{7}{*}{
$F_{n+1}^3+F_n^3-F_{n-1}^3=
\textcolor{blue}{
F_{3n}
}
$
} &
\multirow{7}{*}{$n\geq1$	} &
Lucas \\
&&& \hspace{1cm} \emph{1919 \hspace{2mm} Dickson} \\
&&& \hspace{1cm} \emph{1951 \hspace{2mm} Vorob'ev} \\
&&& \hspace{1cm} \emph{1953 \hspace{2mm} Ginsburg} \\
&&& \hspace{1cm} \emph{1953 \hspace{2mm} Subba Rao} \\
&&& \hspace{1cm} \emph{1972 \hspace{2mm} Brousseau} \\
&&& \hspace{1cm} \emph{2003 \hspace{2mm} Benjamin \& Quinn} \\
\hline
\multirow{3}{*}{1877} &
\multirow{3}{*}{
$5F_n^3=
\textcolor{blue}{
F_{3n}
}
-3(-1)^n
\textcolor{blue}{
F_n
}
$
} &
\multirow{3}{*}{$n\geq0$} &
Lucas \\
&&& \hspace{1cm} \emph{1965 \hspace{2mm} Halton} \\
&&& \hspace{1cm} \emph{2008 \hspace{2mm} Melham} \\
\hline
\multirow{3}{*}{1951} &
\multirow{3}{*}{
$10\sum_{k=0}^n F_k^3=
\textcolor{blue}{
F_{3n+2}
}
+(-1)^{n-1} 6
\textcolor{blue}{
F_{n-1}
}
+5$
} &
\multirow{3}{*}{$n\geq1$} &
Vorob'ev \\
&&& \hspace{1cm} \emph{1993 \hspace{2mm} Clary \& Hemenway} \\
&&& \hspace{1cm} \emph{2003 \hspace{2mm} Benjamin \& Quinn} \\
\hline
\multirow{14}{*}{1953} &
\multirow{3}{*}{
$F_{n+2}^3-3F_n^3+F_{n-2}^3=3
\textcolor{blue}{
F_{3n}
}
$
} &
\multirow{3}{*}{$n\geq2$} &
Ginsburg \\
&&& \hspace{1cm} \emph{1965 \hspace{2mm} Halton} \\
&&& \hspace{1cm} \emph{1972 \hspace{2mm} Brousseau} \\
\cline{2-4}
&
$F_{n+1}^3=F_n F_{n+1} F_{n+2}+(-1)^n
\textcolor{blue}{
F_{n+1}
}
$
&
$n\geq0$ &
\multirow{7}{*}{Subba Rao} \\
\cline{2-3}
&
$4\sum_{k=1}^n F_{2k-1} F_{2k+1} F_{2k+3}=F_{2n+2}^3+7
\textcolor{blue}{
F_{2n+2}
}
-8$
&
\multirow{6}{*}{$n\geq1$} & \\
\cline{2-2}
&
$4\sum_{k=1}^n F_{2k-2} F_{2k} F_{2k+2}=F_{2n+1}^3-7
\textcolor{blue}{
F_{2n+1}
}
+6$
&& \\
\cline{2-2}
&
$4\sum_{k=1}^n F_{2k-1} F_{2k} F_{2k+1}=F_{2n+1}^3+
\textcolor{blue}{
F_{2n+1}
}
-2$
&& \\
\cline{2-2}
&
$4\sum_{k=1}^n F_{2k} F_{2k+1} F_{2k+2}=F_{2n+2}^3-
\textcolor{blue}{
F_{2n+2}
}
$
&& \\
\cline{2-2}
&
$4\sum_{k=1}^n F_{2k}^3=F_{2n+1}^3-3
\textcolor{blue}{
F_{2n+1}
}
+2$
&& \\
\cline{2-2}
&
$10\sum_{k=1}^n F_k F_{k+2} F_{k+4}=F_{3n+8}+(-1)^{n+1} 16
\textcolor{blue}{
F_{n+1}
}
-5$
&& \\
\cline{2-4}
& \multirow{2}{*}{
$4\sum_{k=1}^n F_{k+1}^3=F_{n+2}^3+F_{n+1}^3+(-1)^n 3
\textcolor{blue}{
F_n
}
-2$
} &
\multirow{2}{*}{$n\geq1$} &
Subba Rao \\
&&& \hspace{1cm} \emph{1965 \hspace{2mm} Horadam} \\
\cline{2-4}
& \multirow{2}{*}{
$4\sum_{k=0}^n F_{2k+1}^3=F_{2n+2}^3+3
\textcolor{blue}{
F_{2n+2}
}
$
} &
\multirow{2}{*}{$n\geq0$} &
Subba Rao \\
&&& \hspace{1cm} \emph{1993 \hspace{2mm} Clary \& Hemenway} \\
\hline
1963 &
$2\sum_{k=0}^n(-1)^k F_{k+1}^3 =
\textcolor{blue}{
F_{n+4}
}
+(-1)^n F_n F_{n+2}^2-1$
& $n\geq0$ & Zeitlin \\
\hline
\multirow{4}{*}{1965} &
\multirow{4}{*}{
$2\sum_{k=0}^n F_{k+1}^3=F_{n+1} F_{n+2}^2+(-1)^n
\textcolor{blue}{
F_n
}
+1$
} &
\multirow{4}{*}{$n\geq0$} &
Harris \\
&&& \hspace{1cm} \emph{1968 \hspace{2mm} Pond} \\
&&& \hspace{1cm} \emph{1972 \hspace{2mm} Brousseau} \\
&&& \hspace{1cm} \emph{2004 \hspace{2mm} Benjamin, Carnes \& Cloitre} \\
\hline
1970 &
$F_{n+2}^3-F_{n-1}^3-3F_n F_{n+1} F_{n+2}=
\textcolor{blue}{
F_{3n}
}
$
& $n\geq1$ & Hillman (Padilla) \\
\hline
\multirow{2}{*}{1977} &
$F_n F_{n+3}^2-F_{n+2}^3=(-1)^{n+1}
\textcolor{blue}{
F_{n+1}
}
$
&
\multirow{2}{*}{$n\geq0$} &
\multirow{2}{*}{Hoggatt \& Bergum} \\
\cline{2-2}
&
$F_{n+3} F_n^2-F_{n+1}^3=(-1)^{n+1}
\textcolor{blue}{
F_{n+2}
}
$
&& \\
\hline
\multirow{2}{*}{1993} &
$10\sum_{k=1}^n F_k^3=F_{3n+2}-(-1)^n 6
\textcolor{blue}{
F_{n-1}
}
+5$
&
$n\geq1$ &
\multirow{2}{*}{Clary \& Hemenway} \\
\cline{2-3}
&
$4\sum_{k=1}^n F_{2k}^3=(F_{2n+1}-1)^2 (F_{2n+1}+2)$
&
$n\geq1$ & \\
\hline
2003 &
$F_{n+1} F_{n+2} F_{n+6}-F_{n+3}^3=(-1)^n
\textcolor{blue}{
F_n
}
$
& $n\geq0$ & Melham \\
\hline
\multirow{3}{*}{2005}
&
$F_n F_{n+4} F_{n+5}-F_{n+3}^3=(-1)^{n+1}
\textcolor{blue}{
F_{n+6}
}
$
&
$n\geq0$ &
\multirow{3}{*}{Fairgrieve \& Gould} \\
\cline{2-3}
&
$F_{n-2} F_{n+1}^2-F_n^3=(-1)^{n-1}
\textcolor{blue}{
F_{n-1}
}
$
&
$n\geq2$ & \\
\cline{2-3}
&
$F_{n+2} F_{n-1}^2-F_n^3=(-1)^n
\textcolor{blue}{
F_{n+1}
}
$
&
$n\geq1$ & \\
\hline
2008 &
$8F_n^3=3F_n F_{n-2}^2+3F_{2n} F_{n-2}+
\textcolor{blue}{
F_{3n}
}
$
& $n\geq3$ & Benjamin, Eustis \& Plott \\
\hline
2011 &
$F_{n-3} F_{n+1}^2-F_{n-2}^2 F_{n+3}=4(-1)^n
\textcolor{blue}{
F_n
}
$
& $n\geq3$ & Melham \\
\hline
\multirow{2}{*}{2012} &
\hspace{2mm}
$5F_{n+1}^3-(-1)^n 3
\textcolor{blue}{
F_{n+1}
}
=F_{n+3}^3-F_n^3-3F_{n+1} F_{n+2} F_{n+3}=F_{n+1}^3+F_{n+2}^3-F_n^3$
\hspace{2mm}
&
\multirow{2}{*}{$n\geq0$} &
\multirow{2}{*}{Azarian} \\
\cline{2-2}
&
$3(
\textcolor{blue}{
F_{3n}
}
+4
\textcolor{blue}{
F_{3n+3}
}
)=F_{n+4}^3+F_n^3-3F_{n+2}^3$
&& \\
\hline
2013 &
$
\textcolor{blue}{
F_{3n+3}
}
=F_n F_{n+3}^2+F_{n-1} F_{n+2}^2-F_{n-2} F_{n+1}^2$
& $n\geq2$ & Lang \& Lang \\
\hline
\multirow{3}{*}{2018} &
$F_n F_{n+1}^2-F_n^2 F_{n+2}=(-1)^n
\textcolor{blue}{
F_n
}
$
&
$n\geq0$ &
\multirow{3}{*}{Khomovsky} \\
\cline{2-3}
&
$16
\textcolor{blue}{
F_{3n}
}
=F_{n+3}^3+4F_n^3-F_{n-3}^3$
&
$n\geq3$ & \\
\cline{2-3}
&
$63
\textcolor{blue}{
F_{3n}
}
=F_{n+4}^3-7F_n^3+F_{n-4}^3$
&
$n\geq4$ & \\
\hline
\multirow{4}{*}{} &
$
\textcolor{blue}{
F_{3n}
}
=2F_n^3+3F_{n-1} F_n F_{n+1}$
&
$n\geq1$ & \\
\cline{2-3}
&
$
\textcolor{blue}{
F_{3n+1}
}
=F_{n+1}^3+3F_{n+1} F_n^2-F_n^3$
& $n\geq0$ &
% The task of breaking the url over two lines, inside the tabular environment with the columns spanning multiple rows defeated me. Debate rages in the fora, as to the best way to do this. All the ways seem very complicated, except for using a url shortening service, such as tinyurl, so that the line need not be broken. However, the consensus was that you should not use url shortening services in an academic paper, for various good reasons. Instead I decided to use the \href command instead of the \url command, because this allows you to replace the url with arbitrary text, while retaining the real url as the clickable link. This means that while clicking on the link will still work, copy and pasting the link literally will no longer work, nor will typing it in from a printout. To alleviate this downside, I chose a copy-and-pastable url of the website's homepage, for the link text. Another issue that arose was that the style of the text displayed by \href differs from that of the text displayed by \url. Therefore, for consistency, I used the command \nolinkurl inside the text field of the \href command, so that the style of links would be consistent throughout the article.
\href{http://en.wikipedia.org/wiki/Fibonacci_number\#Other_identities}{\nolinkurl{http://en.wikipedia.org}}
\\
\cline{2-3}
&
$
\textcolor{blue}{
F_{3n+2}
}
=F_{n+1}^3+3F_n F_{n+1}^2+F_n^3$
& $n\geq0$ &
\href{http://mathworld.wolfram.com/FibonacciNumber.html}{\nolinkurl{http://mathworld.wolfram.com}}
\\
\cline{2-3}
&
$
\textcolor{blue}{
F_{3n+3}
}
=2F_{n+1}^3+3F_n F_{n+1}^2+3F_n^2 F_{n+1}$
& $n\geq0$ & \\
\hline
\end{tabular}
}
\end{table}

\section{Conclusions}

%This survey was triggered by a search for cubic Fibonacci identities, as these are the ones which could tile 3D space with cubes/cuboids having sides equal to Fibonacci numbers. Somehow unexpectedly, our investigations have revealed only a rather small number of cubic Fibonacci identities, representing a small fraction of all published Fibonacci identities (most of which are linear or quadratic). Finally, out of these, only a handful of cubic Fibonacci identities are homogeneous. These results have revealed that cubic Fibonacci identities are not so common, not very well-known either, and that discovering new ones has not been a priority.
%
%Our future plans are to rely on the homogeneous cubic Fibonacci identities from Table I (or on similar ones) for constructing 3D Fibonacci spirals. Hence, although the search for 3D Fibonacci spirals is still in its infancy (for a preliminary 3D Fibonacci spiral see \cite{Nagy_Cowell_and_Beiu_2018}), quite a few alternatives are clearly feasible.

While Fibonacci numbers appear in nature \cite{Mitchison_1977}, \cite{Livio_2002}, they have also been shown to optimize certain man-made systems. As an example we have designed ourselves \cite{Beiu_1999}, a class of adder circuits implemented using threshold logic gates was shown to be VLSI-optimal as the minimal weights (of the threshold logic gates) are the Fibonacci numbers, which minimizes area while also maximizing speed. 

This survey was triggered by a search for cubic Fibonacci identities, as these are the ones that might tile 3D space with cubes/cuboids having sides equal to Fibonacci numbers. Somewhat unexpectedly, our investigations have revealed only a rather small number of cubic Fibonacci identities, representing a small fraction of all published Fibonacci identities (most of which are linear or quadratic). Finally, out of these, only a handful of cubic Fibonacci identities are homogeneous. These results have revealed that cubic Fibonacci identities are not so common, not very well-known either, and that discovering new ones has not been a priority.

\begin{figure}
\begin{subfigure}[b]{0.19\textwidth}
\includegraphics[width=\linewidth]{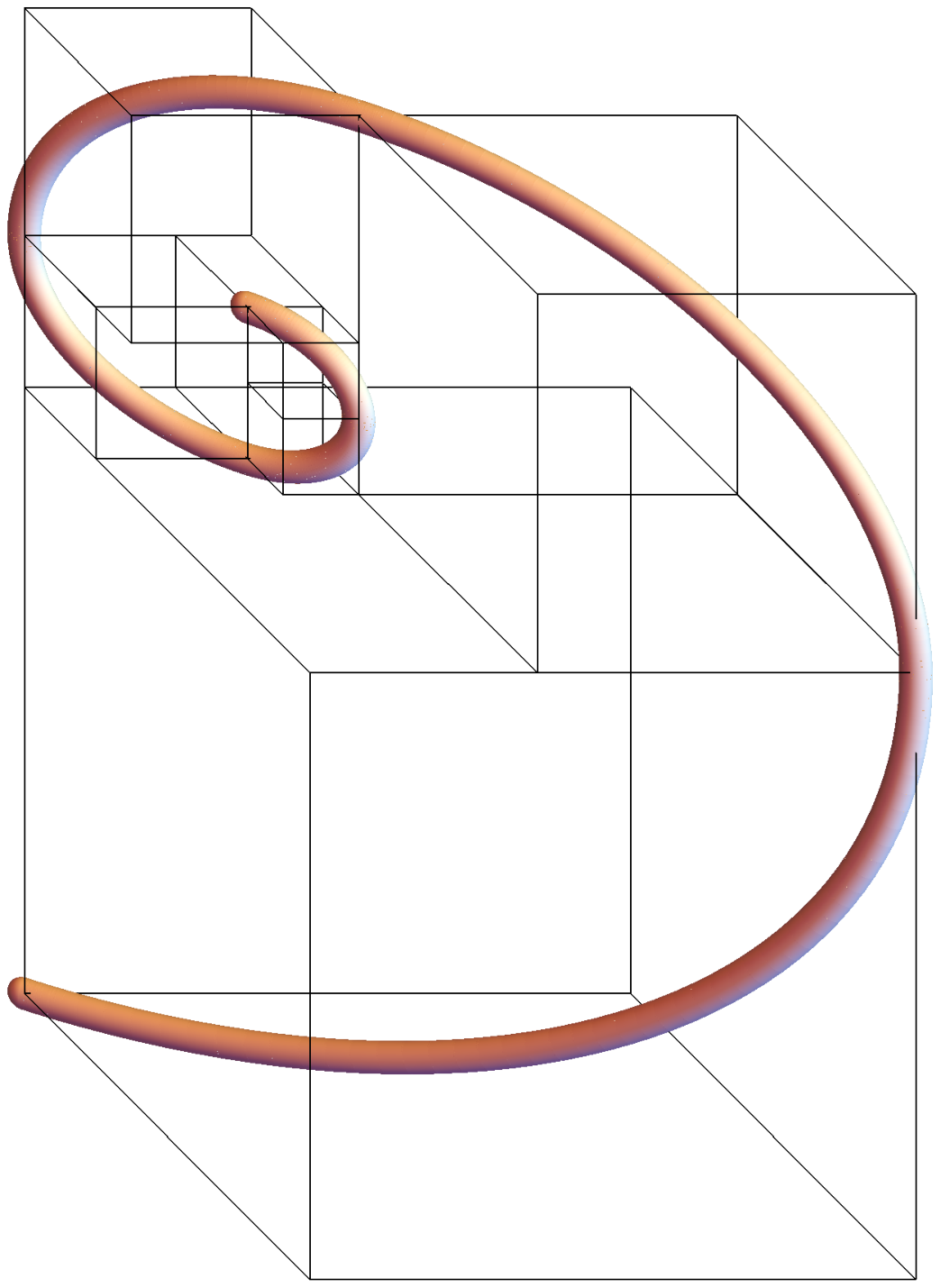}
\caption{}
\label{fig: VBplanarfibtubespiral}
\end{subfigure}
\hspace{\fill}
\begin{subfigure}[b]{0.19\textwidth}
\includegraphics[width=\linewidth]{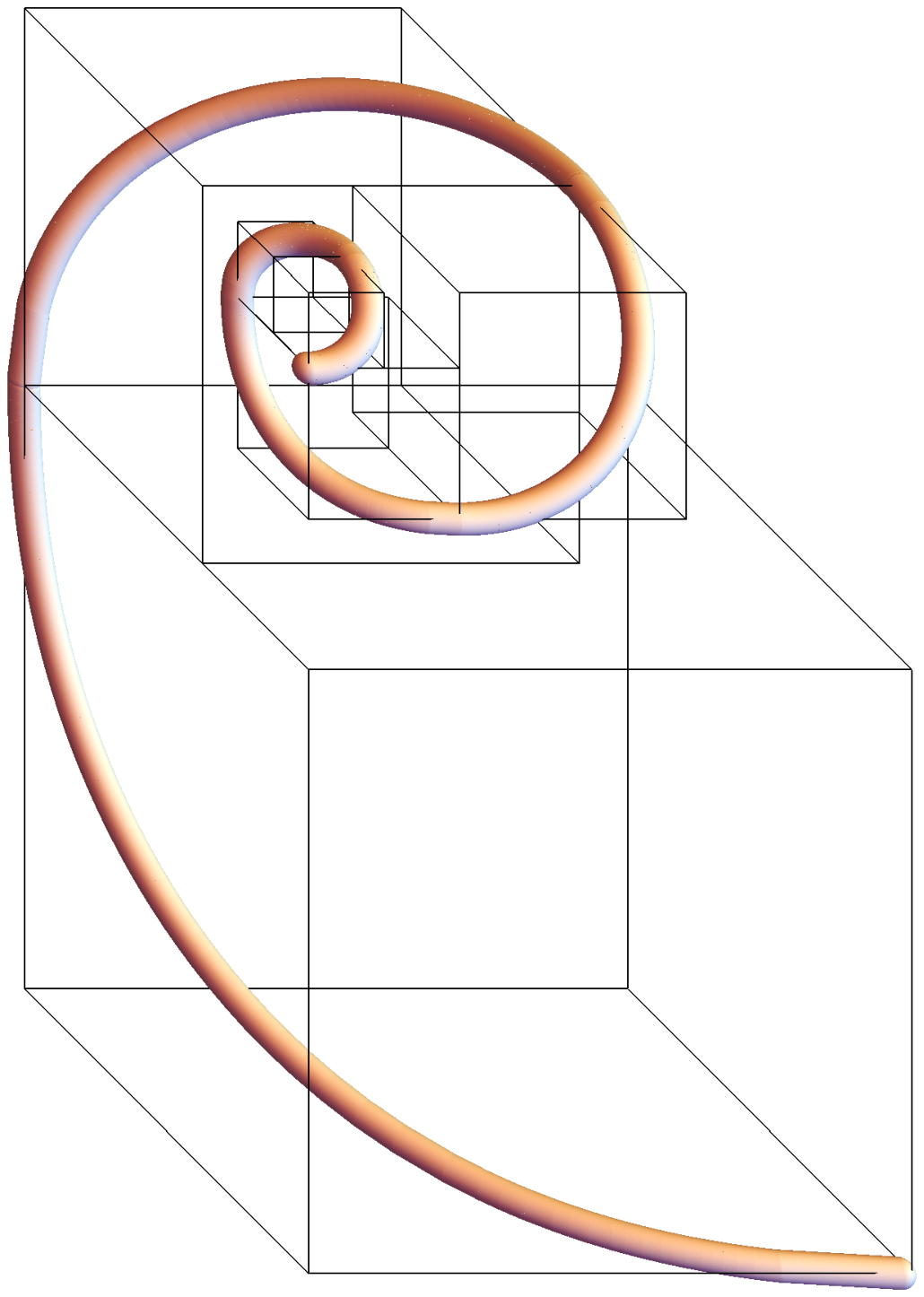}
\caption{}
\label{fig: diagonalpadovanstylefibtubespiralv2}
\end{subfigure}
\hspace{\fill}
\begin{subfigure}[b]{0.19\textwidth}
\includegraphics[width=\linewidth]{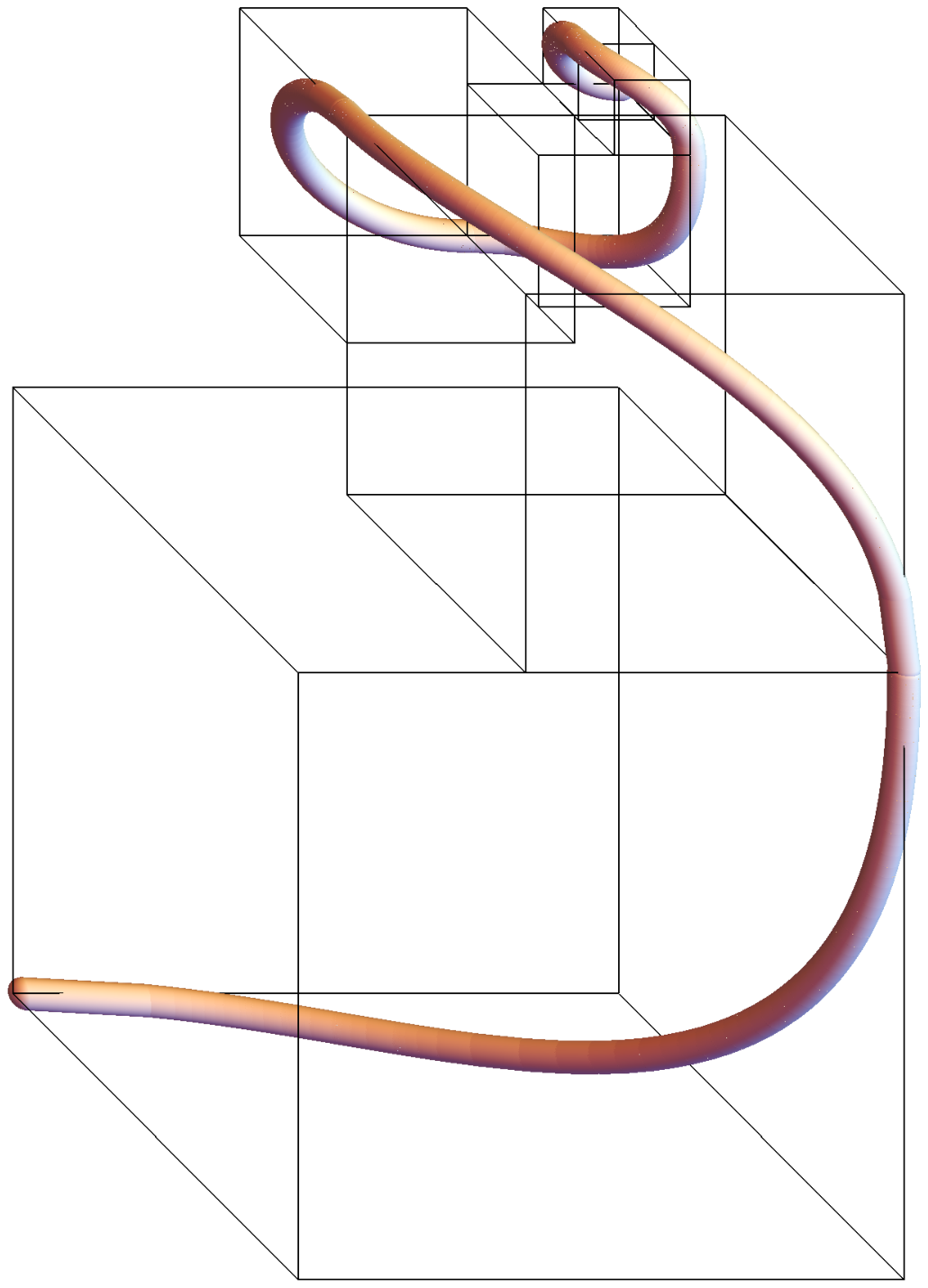}
\caption{}
\label{fig: diagonalpadovanstylefibtubespiralv4}
\end{subfigure}
\hspace{\fill}
\begin{subfigure}[b]{0.19\textwidth}
\includegraphics[width=\linewidth]{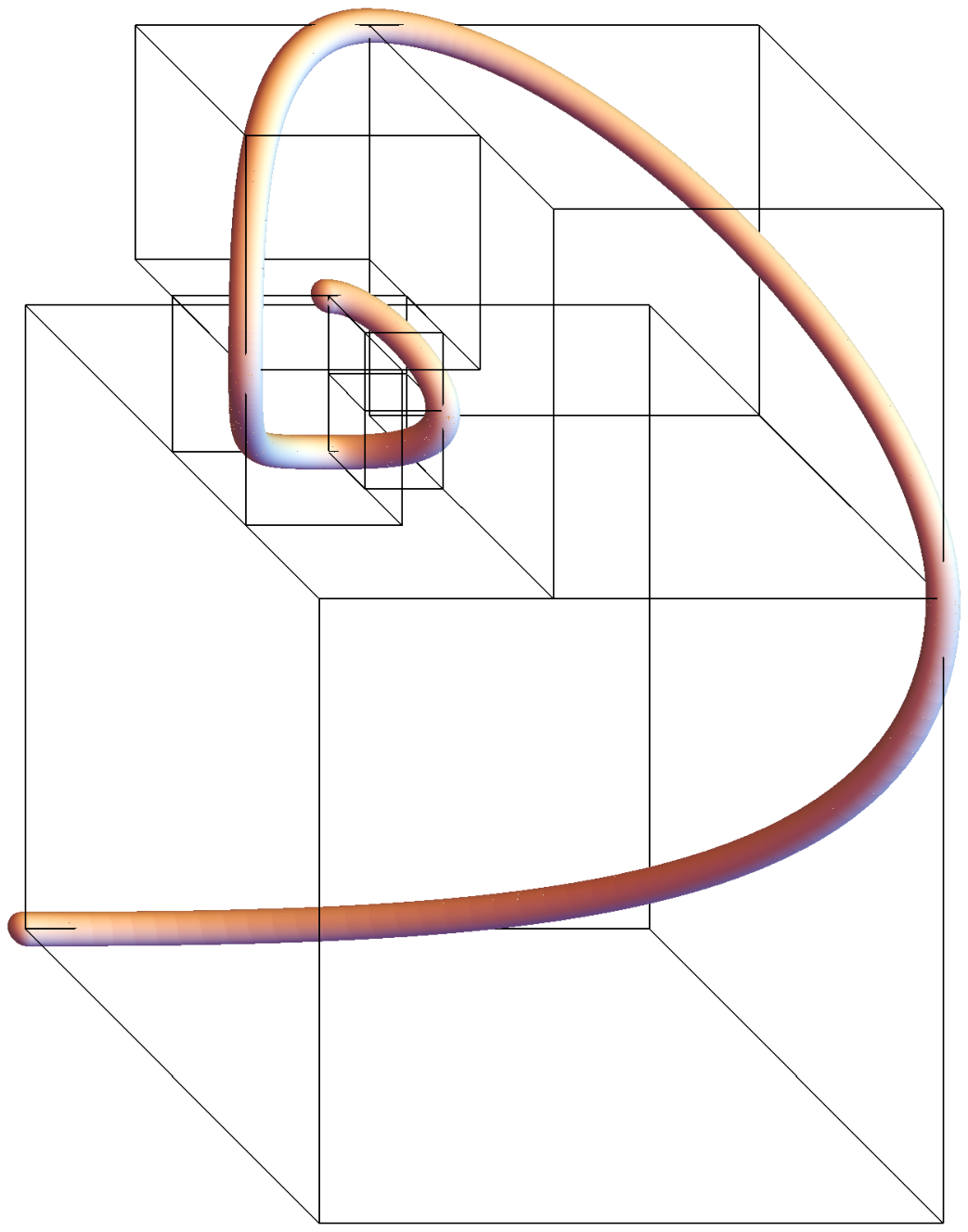}
\caption{}
\label{fig: VBfibtubespiralv1}
\end{subfigure}
\hspace{\fill}
\begin{subfigure}[b]{0.19\textwidth}
\includegraphics[width=\linewidth]{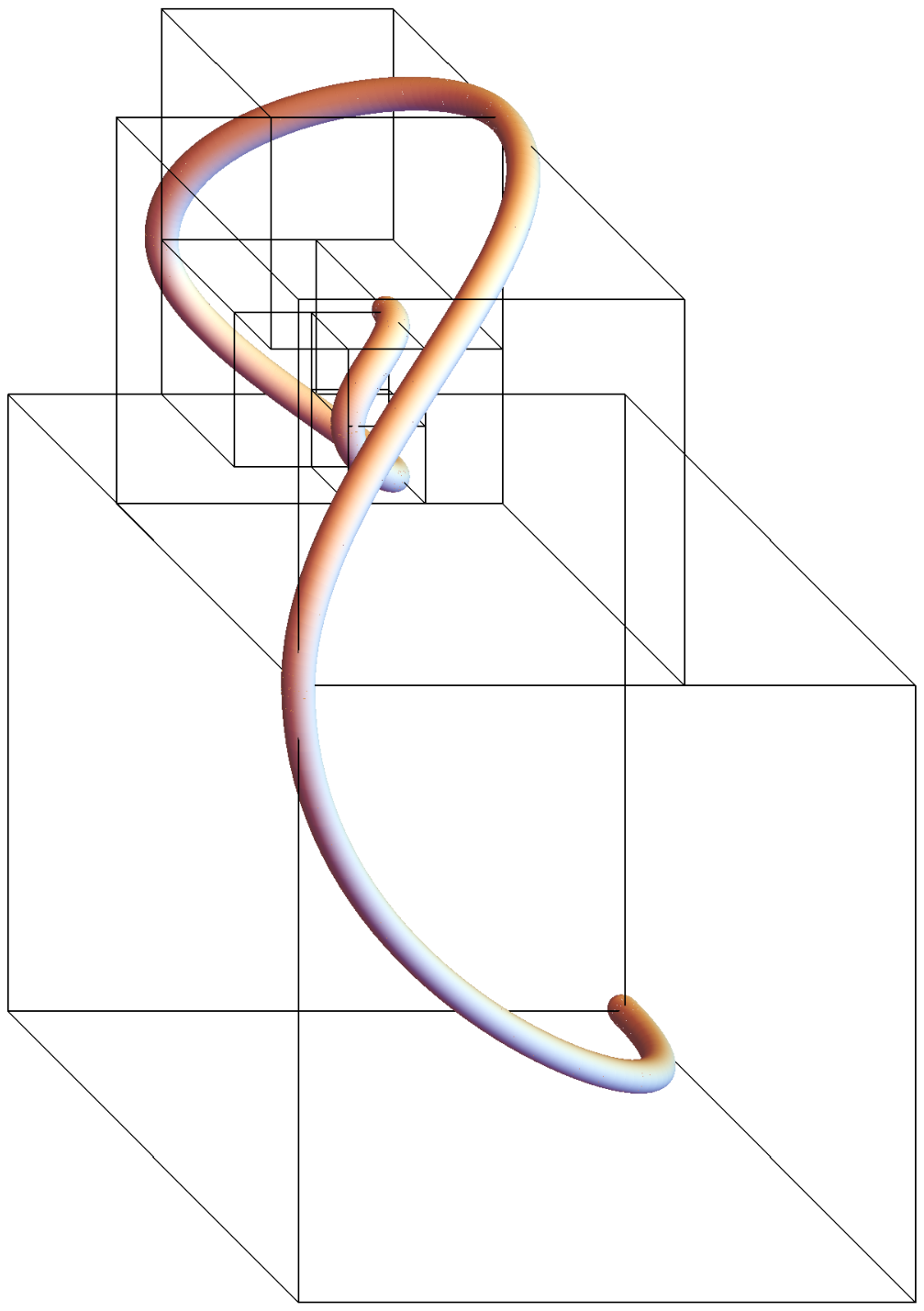}
\caption{}
\label{fig: VBfibtubespiralv4}
\end{subfigure}
%\hspace{\fill}
\caption{(a) 3D Fibonacci spiral from \cite{Nagy_Cowell_and_Beiu_2018} (reproduced with permission) based on $(M^{(1)}_n,b^{(1)}_n)$; (b)-(c) Two views of a spiral based on $(M^{(2)}_n,b^{(2)}_n)$; (d)-(e) Two views of a spiral based on $(M^{(3)}_n,b^{(3)}_n)$.}
\label{fig: spirals}
\end{figure}

Our short term plans are to construct 3D Fibonacci spirals relying on the homogeneous cubic Fibonacci identities from Table \ref{tab: I}. These homogeneous cubic identities suggest quite a large number of different tilings of 3D space using Fibonacci cuboids -- especially when distinct such identities are used in combination. It also means that many 3D Fibonacci spirals are possible. A first 3D Fibonacci spiral based on eq.~\eqref{eq: 31} (see \cite{Pond_1968} as well as \cite{Benjamin_Carnes_and_Cloitre_2009}, \cite{Lang_2012}) was presented for the first time in \cite{Nagy_Cowell_and_Beiu_2018}, and can be seen in Fig.~\ref{fig: VBplanarfibtubespiral}.~This is a 3D Fibonacci spiral analogous to the 3D Padovan one (see \cite{OEIS_A000931}, \cite{Padovan_1994}), but having a growth factor of $\varphi \sqrt{2} \approx 2.288245$, where $\varphi = (1+\sqrt{5})/2$ is the golden ratio. Unfortunately, exactly like the 3D Padovan spiral, which is in fact in 2D \cite{Stewart_1996}, this 3D Fibonacci spiral is also in 2D (on a diagonal plane). Still, ``proper'' 3D Fibonacci spirals can be constructed as suggested by the 3D packing shown in Figs.~\ref{fig: diagonalpadovanstylefibtubespiralv2} and \ref{fig: diagonalpadovanstylefibtubespiralv4}, and by the one shown in Figs.~\ref{fig: VBfibtubespiralv1} and \ref{fig: VBfibtubespiralv4}.~Projecting this last spiral orthogonally, one recovers exactly the 2D Fibonacci spiral, i.e. with growth factor $\varphi$. The arrangements of Fibonacci cubes in Fig.~\ref{fig: spirals} are quasi self-similar fractals defined by a sequence of affine maps $C_{n+1} = M^{(i)}_n C_n + b^{(i)}_n$, $n\geq1$, where $C_1=[0,1]^3$, the unit cube.
Precisely, for all $n \geq 1$,
\begin{align*}
M^{(1)}_n
&=
\frac{F_{n+1}}{F_n}
\left (
\begin{array}{ccc}
(-1)^{n+1}		&	0					&	0 \\
0					&	(-1)^{n+1}		&	0 \\
0					&	0					&	(-1)^n
\end{array}
\right ),
%\quad
&
b^{(1)}_n
&=
\left \{
\begin{array}{cc}
(0,0,0)^T,
&
n \text{ odd}, \\
\left (
1+\frac{F_{n+1}}{F_n}, 1+\frac{F_{n+1}}{F_n}, 1-\frac{F_{n+1}}{F_n}
\right)^T,
& n \text{ even},
\end{array}
\right .
\\
M^{(2)}_n
&=
\frac{F_{n+1}}{F_n}
\left (
\begin{array}{ccc}
0	&	-1	&	0 \\
0	&	0	&	1 \\
-1	&	0	&	0
\end{array}
\right ),
&
b^{(2)}_n
&=
\frac{F_{n+1}}{F_n}
\left (
\begin{array}{c}
1 \\
-1 \\
1
\end{array}
\right ),
%(1,-1,1)^T,
\\
M^{(3)}_n
&=
\frac{F_{n+1}}{F_n}
\left (
\begin{array}{ccc}
1					&	0					&	0 \\
0					&	(-1)^{n+1}		&	0 \\
0					&	0					&	(-1)^n
\end{array}
\right ),
&
b^{(3)}_n
&=
\left \{
\begin{array}{cc}
(0,0,0)^T,
&
n \text{ odd}, \\
\left (
1-\frac{F_{n+1}}{F_n}, 1+\frac{F_{n+1}}{F_n}, 1-\frac{F_{n+1}}{F_n}
\right)^T,
& n \text{ even}.
\end{array}
\right .
\end{align*}

Our long term plans are to use such identities for 3D packing (a direction of research started by \cite{Kepler_1611}). We expect that these could be useful for explaining particular behaviors of (complex) ionic solutions in confined spaces \cite{Munteanu_Cavaco-Paulo_and_Beiu_2017} - as atomic radii are related to Fibonacci numbers through the 3D packing of protons and neutrons \cite{Kepler_1611}, with covalent radii following suit. Additionally, such cubic identities also lead to fractal structures \cite{Mandelbrot_1982} (and fractional calculus), while combining them with physical/chemical constraints might help explain certain protein structures \cite{Enright_and_Leitner_2015}, and the fractality of cortical actin meshwork \cite{Sadegh_Higgins_Mannion_Tamkun_and_Krapf_2017}. That is why we believe that, at least the homogenous cubic Fibonacci identities reported in Table \ref{tab: I} as well as the homogenous ones from Table \ref{tab: III}, could prove useful for future research related to several different fields.

\begin{table}%[!h]
%\label{tab: III}
\centering
\caption{Selection of General Fibonacci Identities (In Chronological Order)}
\label{tab: III}
\resizebox{\textwidth}{!}{
\begin{tabular}{|c|c|c|l|}
\hline
\hline
Year & Identity & Condition & Author(s) \\
\hline
\hline
1963 &
$F_{r+1} F_{s+1} F_{t+1}+F_r F_s F_t-F_{r-1} F_{s-1} F_{t-1}=
\textcolor{blue}{
F_{r+s+t}
}
$
&
$r,s,t\geq1$ &
Hoggatt \\
\hline
\multirow{2}{*}{
1964
} &
$F_{n+a} F_{n+b} F_{n+c}-F_n F_{n+a} F_{n+b+c}+F_n F_{n+b} F_{n+c+a}-F_n F_{n+c} F_{n+a+b}$
&
\multirow{2}{*}{
$n,a,b,c\geq0$
} &
\multirow{2}{*}{Gould} \\
&
$=(-1)^n (F_a F_b F_{n+c}-F_c F_a F_{n+b}+F_b F_c F_{n+a})$
&& \\
\hline
1965 &
$F_m^k F_n=(-1)^{kr} \sum_{h=0}^k C_k^h (-1)^h F_r^h F_{r+m}^{k-h} F_{n+kr+hm}]$
&
$k,m,n,r\geq0$ &
Halton \\
\hline
2006 &
\hspace{2mm}
$F_{l+1} F_{m+1} F_{n+1}=3F_l F_m F_n+6F_{l-1} F_{m-1} F_{n-1}-3F_{l-2} F_{m-2} F_{n-2}-F_{l-3} F_{m-3} F_{n-3}$
\hspace{2mm}
&
$l,m,n\geq3$ &
Ollerton \\
\hline
\multirow{3}{*}{2008} &
\multirow{3}{*}{
$\sum_{k=0}^n C_n^k F_{qk+p} F_{q-2}^{n-k}=F_q^n F_{2n+p}$
} &
$n\geq1$ &
\multirow{3}{*}{Benjamin, Eustis \& Plott} \\
&& $p\geq0$ & \\
&& $q\geq3$ & \\
\hline
2012 &
$F_{k,n+4}^3=(k^3+2k) F_{k,n+3}^3+(k^4+3k^2+2) F_{k,n+2}^3-(k^3+2k) F_{k,n+1}^3-F_{k,n}^3$
&
$n\geq0$ &
Yazlik, Yilmaz \& Taskara \\
\hline
\multirow{2}{*}{2018} &
\multirow{2}{*}{
$F_{2l} (F_{n+m}^2-F_{n-m}^2 )=F_{2m} (F_{n+l}^2-F_{n-l}^2)$
} &
$l,m\geq1$ &
\multirow{2}{*}{Khomovsky} \\
&& $n \geq \max(l,m)$ & \\
\hline
\multirow{3}{*}
&
\multirow{3}{*}{
$
\textcolor{blue}{
F_{kn+m}
}
=\sum_{i=0}^k C_k^i F_{m-i} F_n^i F_{n+1}^{k-i}$
} &
$k\geq0$ &
%\url{http://mathworld.wolfram.com/FibonacciNumber.html} \\
\href{http://mathworld.wolfram.com/FibonacciNumber.html}{\nolinkurl{http://mathworld.wolfram.com}} \\
&& $n\geq1$ &
\hspace{1cm} \emph{2003 \hspace{2mm} Mihailov} \\
&& $m\geq k+1$ &
\hspace{1cm} \emph{2018 \hspace{2mm} Nagy, Cowell \& Beiu} \\
\hline
\end{tabular}
}
\end{table}

%\vspace*{5cm}
%\par
%\ \\
%\vspace{5cm}
% I've been trying for quite a long time, and so far have failed, to make a space between the end of table 3, and the start of the figures.

% The following code is taken from Stackexchange, here: https://tex.stackexchange.com/a/132604/176433

% This \clearpage command seems to be necessary, to force the bibliography to go at the end of the document. Otherwise Table II is put after the bibliography. Or if you remove the [!h] position options for the table environments, as I now have, then tables II and III appear part of the way through the bibliography ! LaTeX is notorious for its ``intelligent'' floating rules...

\vspace{1cm}

\noindent
Department of Mathematics and Computer Science, ``Aurel Vlaicu" University of Arad, Str. Elena Dr\u{a}goi 2-4, RO-310330 Arad, Romania. \\
mariana.nagy@uav.ro, simon.cowell@uav.ro, valeriu.beiu@uav.ro

\vspace{1cm}

\textbf{Mariana Nagy} received the B.Sc. degree in computer engineering and Ph.D. in engineering, from ``Politehnica'' University of Timisoara, Timisoara, Romania in 1985, and respectively 1996, and the M.S. degree in economics from ``Aurel Vlaicu'' University of Arad, Arad, Romania in 2006.

From 1985 to 1992 she was with ICSIT-MU (SIMTEX) doing research on machine-tools, from 1992 to 1998 she was with SC ``Sanevit'' SA, and since 1999 she is with the ``Aurel Vlaicu'' University of Arad, Arad, Romania, in various positions. At present she is professor with the Department of Mathematics and Computer Science and the Dean of the Faculty of Exact Sciences. She published 8 books and over 70 papers, and was PI/co-PI on 16 EU projects. Her research interests include mathematical modeling and computer applications mainly for economic phenomena and medical statistics. 

\vspace{5mm}

\textbf{Simon R. Cowell} received the B.Sc. and Ph.D. degrees in mathematics from the University of Wales, Swansea, UK in 1999, and respectively 2003, and the Ph.D. degree in mathematics from the University of Missouri at Columbia, USA in 2009.

From 2009 to 2017, he was an Assistant Professor with the United Arab Emirates University (UAE). Since 2017, he has been a Senior Researcher/Lecturer with the Department of Mathematics and Computer Science, ``Aurel Vlaicu'' University of Arad, Arad, Romania. His research interests include reliability, mathematics for medical imaging, complex analysis and Banach space theory.

He is a member of the London Mathematical Society. 

\vspace{5mm}

\textbf{Valeriu Beiu} received the M.Sc. degree in computer engineering from the University ``Politehnica'' of Bucharest (UPB), Bucharest, Romania, in 1980, and the Ph.D. degree (summa cum laude) in electrical engineering from the Katholieke Universiteit Leuven (KUL), Leuven, Belgium, in 1994.

He was with the Research Institute for Computer Techniques, Romania, working on high-speed CPUs/FPUs, prior to returning to UPB, while later he was with KUL (1991-1994), King’s College London (1994-1996), Los Alamos National Laboratory (LANL, 1996-1998), RN2R/Rose Research (1998-2001), Washington State University (2001-2005), and the United Arab Emirates University (2005-2015), where he was the Associate Dean for Research and Graduate Studies (2006-2011). Since 2015, he joined the ``Aurel Vlaicu'' University of Arad (UAV), Arad, Romania, where he is currently Professor and Director. He was a PI/co-PI of research grants totaling over US\$ 51M, gave over 200 invited talks, organized over 120 conferences, chaired over 60 sessions, has edited 2 books and authored over 240 journal/conference articles (30 invited), 8 chapters, and holds 11 patents. His research interests focus on biological-inspired nano-circuits and brain-inspired nano-architectures for VLSI-efficient designs (ultra-low power and highly reliable).

Dr. Beiu received five Fellowships, including Fulbright (1991), EU Human Capital \& Mobility (1994-1996), LANL Director's Postdoc (1996-1998), and Fellow of Rose Research (1999-2001). He is a senior member of IEEE (since 1996) and a member of ACM, INNS, ENNS, MCFA, IEEE Computer Science Task Force on Nanoarchitectures (2005), IEEE Emerging Technology Group on Nanoscale Communications (2009), and was a member of the SRC-NNI Working Group on Novel Nano-Architectures (2003-2006). He was the recipient of seven Best Paper Awards, was the Program Chairman of the IEEE LANL Section (1997), and has served as an Associate Editor of the IEEE Transactions on Neural Networks (2005-2008), Nano Communication Networks (2010-2015), and the IEEE Transactions on Very Large Scale Integration Systems (2011-2015).

\vspace{1cm}

\end{document}